\documentclass[a4paper,twoside]{article}
\usepackage{amssymb}
\usepackage[latin1]{inputenc}
\usepackage[T1]{fontenc}
\usepackage[final]{lpretex}
\usepackage[all]{xy}
\usepackage{amscd}
\usepackage{title}

\DocVersion[Sun Sep 22 20:23:58 CEST 2002]{1}

\newcommand\Ab[3][]{{\Cal H}^{#3}_{#2\if\relax#1\relax\else,#1\fi}}

\newcommand\tAb[2]{{\Cal H}^{t,#2}_{#1}}
\newcommand\sAb[2]{{\Cal H}^{s,#2}_{#1}}
\newcommand\uAb[2]{{\Cal U}^{#2}_{#1}}
\newcommand\vAb[2]{{\Cal V}^{#2}_{#1}}
\newcommand\wAb[2]{{\Cal W}^{#2}_{#1}}
\newcommand\Hyp[2]{{\Cal H}_{#1,#2}}
\newcommand\Pol[1]{{\Cal P}(#1)}
\newcommand\Aff{\symb{Aff}}
\newcommand\Conf[2]{\symb{Conf}^{#2}({#1})}
\newcommand\oConf[2]{\symb{\widetilde{Conf}}^{#2}({#1})}

\begin{document}

\title{On Abel's hyperelliptic curves}

\author{T. Ekedahl}
\address{Department of Mathematics\\
 Stockholm University\\
 S-106 91  Stockholm\\
Sweden}

\maketitle

\begin{dedication}
Dedicated to Niels Henrik Abel on his 200'th birthday.
\end{dedication}

\begin{abstract}
In this note we discuss a class of hyperelliptic curves introduced by Abel (cf.,
\cite{abel26::sur+r+r}). After some indications of the context in which he
introduced them and a description of his main result we give some results on the
moduli space of such curves.

In particular we compute the dimension of it at each of its points as well as
giving a combinatorial formula for the number of components.
\end{abstract}

%;;%%Introduction
In the paper \cite{abel26::sur+r+r} takes up a very special case of the problem
of deciding when a rational differential form is the logarithmic differential of
a rational function. Even though it is easy to imagine that the problem later
led Abel to his famous Paris dissertation the solution that Abel proposed is
very special and quite different from the approach he chose later. Nevertheless,
the present article is concerned with investigating the problem of [loc.\ cit.].

A characteristic, and for the time quite unusual, feature of all of Abel's work
is his insistence on treating general cases rather than special examples. It
thus seems entirely fitting to study the moduli problem of all the solutions to
his proposed problem and we shall indeed consider a number of aspects of the
moduli spaces (or more precisely stacks) that classify his solutions.

We start by introducing a (small) number of variants of an attempt to formulate
Abel's condition as a moduli problem. In particular Abel's solution of his
problem in terms of a polynomial type Pell's equation appears not as the moduli
problem that has been chosen as the central moduli problem of this article but
as a chart for it. We then study the relations between these variants, the end
result being that they are indeed closely related.

After a discussion of the case of genus $0$ and $1$ we give a reformulation of
the moduli problem in terms of families of maps between genus zero curves. We
then proceed to make an infinitesimal study of the moduli problem which allows
us to conclude that it is smooth in characteristic $0$. We then go on to study a
Lyashko-Looijenga type map and show that it is covering map. This allows us to
give a topological covering space type description of the moduli stack which in
particular gives us a combinatorial description of the set of components of the
moduli space.

We shall, except for the last section, adopt a purely algebraic approach. Apart
from reasons of taste there are some arguments in favour of such an
approach. The reader's attention should be particularly directed to Theorem
\ref{Abel deformations} where we shall discover that some naturally defined
``equi-ramification strata'' turn out to be non-reduced. It seems likely that
the multiplicity with which those strata appear is significant.

We shall also use the language of algebraic stacks. This may seem unnecessary
particularly as our stacks are very close to being spaces (cf., \(fixed
points)). However, I claim that it is the technically most convenient as well as
most intuitive way of doing things. In particular when defining maps between
solutions to moduli problems representing these solutions as stacks means that
in order to define maps between them one may often follow the path of first
deciding what the map should do on points and then verify that this pointwise
construction is natural enough so that it makes sense for families of
objects. This is in fact what we shall do most of the time. Sometimes, however,
we shall discover that some choices that were made in the point case can not be
made in the case of a family and we shall then have to incorporate those choices
in the definition of the moduli problem. This will lead for instance to the
three slightly different versions of Abel curve.
\begin{conventions}
By a \Definition{monic polynomial} we shall mean a polynomial in one variable
whose highest degree coefficient is equal to $1$. Such a polynomial will be said
to be \Definition{normalised} if its next to highest degree coefficient is equal
to $0$.

As we shall deal extensively with stacks it seems natural to use the term
`scheme' to denote an algebraic space and hence by `locally' mean `locally in
the \'etale topology'.\footnote{Note that in practice the only difference
between ordinary schemes and algebraic spaces is that for the latter the Zariski
topology is not available.} Though we shall do so this is not strictly necessary
however and then `locally' may at times be interpreted as `locally in the
Zariski topology' though consistently using the \'etale topology will always
work.

We have made only a token attempt at formulating our results in arbitrary
characteristics.\footnote{It will be clear that if the characteristic is large
enough with respect to the degree $n$ the situation will be similar to that of
characteristic $0$.} Starting with section \ref{sec:Abel curves} all our schemes
and stacks will be over $\Sp\Z[1/2]$ and starting with section
\ref{sec:infinitesimal} we shall work exclusively in characteristic zero, this
will also be true at the end of the preliminary section \ref{def:prelims} and at
points in section \ref{sec:Abel curves} (which will be explicitly spelled out).

As usual a \Definition{multiset} is a set whose members are counted with certain
multiplicities, formally it is a set provided with a multiplicity function from
it to the integers $\ge 1$. If $S$ is a multiset we shall use $|S|$ to denote
the domain of the multiplicity function and $\mu_S$ for the multiplicity
function itself. We shall use set-theoretic notation when dealing with
multisets:
\begin{itemize}
\item $S:=\{1,1,1,2,2,3\}$ will denote the multiset for which $|S|=\{1,2,3\}$
and $\mu_S(n)=4-n$.

\item $\sum_{s \in S} s^2$ should be interpreted as $\sum_{s \in
|S|}\mu_S(s)s^2$, i.e., $3\cdot1^2+2\cdot2^2+3^2$.

\item Similarly, \set{\lfloor s/2\rfloor}{s \in S} should be interpreted as
$\{0,0,0,1,1,1\}$.
\end{itemize}
A multiset $S$ is finite if $|S|$ and then its \Definition{cardinality} equals
$\sum_{s \in S}1$. A multiset $S$ is said to be a submultiset of the multiset
$T$ if the $|S| \subseteq |T|$ and $\mu_S(s) \le \mu_T(s)$ for all $s \in |S|$.
\end{conventions}

%;;%%Preliminaries Section def:prelims
\begin{section}{Preliminaries}\label{def:prelims}

We shall sometimes speak about the universal object over a stack which
classifies some type of geometric object. Note that, contrary to the case when
the moduli problem of classifying such objects is representable by a scheme,
this is somewhat ambiguous and is not quite as strong. Firstly, for a family of
objects over $S$, the family may not be the pullback of the universal family but
is so only locally on $S$. Secondly, the universal object is not unique; two
such objects are only locally isomorphic. Thirdly, a universal object may in
fact not even exist over the stack itself but only. It would be more proper to
speak about the stack of universal objects but we shall allow ourselves the
luxury of not doing that. The first phenomena are shown quite clearly in the
case of the classifying stack, $BG$, of a finite group $G$. A universal object
is given by the trivial $G$-torsor and a non-trivial $G$-torsor over $S$ is of
course not the pullback of the trivial one. In fact, any $G$-torsor over the
base is universal and there may very well be non-trivial $G$-torsors over the base.

Assume that $X \to S$ is a scheme and $X \to \P^1\times S$ an $S$-morphism. Let
$C$ and $D$ be the schematic inverse images of $0\times s$ and $\infty\times S$
and assume $C$ and $D$ are Cartier divisors. We shall repeatedly use the
(obvious) fact that such a morphism is the same thing as an isomorphism
$\cO_X(C) \riso \cO_X(D)$.

Let $S$ be a scheme. A line bundle $\cL$ and a trivialisation $\varphi$ of
$\cL^2$ will be called an \Definition{involutive line bundle}. Consider further
$\P(\cO_S\Dsum\cL) \to S$, the projective bundle on the vector bundle
$\cO_S\Dsum\cL$, the two sections $\infty$ and $\mathbf 0$ associated to the two
projections of $\cO_S\Dsum\cL$ and the involution $\sigma$ of $\P(\cO_S\Dsum\cL)
\to S$ defined as the composite of the map $\P(\cO_S\Dsum\cL) \to
\P(\cL\Dsum\cO_S)$ that switches the two factors, the standard identification and
distributivity
$\P(\cL\Dsum\cO_S)=\P((\cL\Dsum\cO_S)\Tensor\cL)=\P(\cL^2\Dsum\cL)$ and
$\varphi$ applied to the first factor $\P(\cL^2\Dsum\cL) \to
\P(\cO_S\Dsum\cL)$. We shall call the data
$(\P(\cO_S\Dsum\cL),\mathbf{0},\infty,\sigma)$ the \Definition{involutive
projective bundle} associated to the involutive line bundle $(\cL,\varphi)$ and
denote it $\P_{\cL}$.
\begin{example}
If $\cL=\cO_S$ and $\varphi=\lambda \in \cO_S^\times$, then the involutive
bundle is $(\P^1\times S,0\times S,\infty\times S,x \mapsto \lambda/x)$. Locally
this is the general situation.
\end{example}
Note also that the fixed point locus of $\sigma$ is a double covering of the
base that is isomorphic to the double cover associated to $\cL$ and $\varphi$;
something which is seen for instance by using the local description just
given. We shall call it the \Definition{involutive locus}.

Similarly to the remark above, an $S$-morphism $X \to \P_{\cL}$ such
that the inverse images of $\infty$ and $\mathbf 0$ are Cartier divisors $C$ and
$D$ is the same thing as an isomorphism $\cO_X(C) \riso \cO_X(D)\Tensor\cL$.

The involutive bundle will be said to be \Definition{split} if one is given a
trivialisation of $\cL$ for which $\varphi$ becomes the identity. Then the
involutive projective bundle is identified with $\P^1\times S$ in such a way
that $0$ corresponds to the zero section, $\infty$ to the section at infinity
and the involutive locus is given by \set{(s\co t)}{s^2=t^2} which when $2$
is invertible is $\{(\pm1\co 1)\}$.

We may explicitly construct the quotient of \map{\tau}{\P(\cO_S\Dsum\cL)}{S} by
the action of $\sigma$ in the following way. We define an $S$-map
$\P(\cO_S\Dsum\cL) \to \P(\cO_S\Dsum\cL)$ by giving $\cO_{\P(\cO_S\Dsum\cL)}(2)$
as a quotient of $\tau^*(\cO_S\Dsum\cL)$. By adjunction giving such a map is the
same as giving a map $\cO_S\Dsum\cL \to \tau_*\cO(1)=S^2(\cO_S\Dsum\cL)$. We do
this by mapping $1$ of the $\cO$-factor to $1\tensor 1\oplus \varphi(1)$ in
$S^2\cO_S\Dsum \cL^{\tensor 2}\subset S^2(\cO_S\Dsum\cL)$ and the $\cL$-factor
to $\cO_S\Tensor\cL \subset S^2(\cO_S\Dsum\cL)$ through $1\tensor id$.  In the
local normal form above -- homogenised -- this map is given by $(x\co y) \mapsto
(x^2+\lambda y^2\co xy)$ which evidently has no base points, i.e., it is
surjective and hence gives a map
\map{\pi}{\P(\cO_S\Dsum\cL)}{\P(\cO_S\Dsum\cL)}. As $\sigma$ locally has the
form $(x\co y) \mapsto (\lambda y\co x)$ it is clear that $\pi$ is equivariant
with trivial action on the target. Using again the local form it is easily
verified that it is the quotient map. We shall speak of it as the
\Definition{involutive quotient map} associated to $\cL$ and $\varphi$. Note
that the involutive locus maps to the image under $x \mapsto 2x$ of
itself. For this reason, starting with section \ref{sec:Abel curves}, we shall
instead use $(x\co y) \mapsto (x^2+\lambda y^2\co 2xy)$ as quotient map so
that the involutive locus is mapped to itself.

Seen from the side of its target the involutive quotient map $\pi$ is a double
covering. Restricting ourselves to the case when $2$ is invertible we may
describe this covering as follows. We get a map $\cO_{\P(\cO_S\Dsum\cL)}(-1) \to
\pi_*\cO_{\P(\cO_S\Dsum\cL)}$: By adjunction it corresponds to a map
$\pi^*\cO_{\P(\cO_S\Dsum\cL)}(-1) \to \cO_{\P(\cO_S\Dsum\cL)}$ and by
construction we have
$\pi^*\cO_{\P(\cO_S\Dsum\cL)}(-1)=\cO_{\P(\cO_S\Dsum\cL)}(-2)$ so that such a
map corresponds to a section of $\cO_{\P(\cO_S\Dsum\cL)}(2)$, i.e., a section of
$S^2(\cO_S\Dsum\cL)$ and we choose $\tensor 1\oplus -\varphi(1)$. In the local
form above the section $\tensor 1\oplus -2\varphi(1)$ corresponds to
$1/2(x^2-\lambda y^2)$. From that it is easily verified that the map
$\cO_{\P(\cO_S\Dsum\cL)}(-1) \to \pi_*\cO_{\P(\cO_S\Dsum\cL)}$ is injective and
has as image the $-1$-eigenspace of $\sigma$. The double cover $\pi$ is now
determined by the square map
$\cO_{\P(\cO_S\Dsum\cL)}(-1)\Tensor\cO_{\P(\cO_S\Dsum\cL)}(-1) \to
\cO_{\P(\cO_S\Dsum\cL)}$, i.e., a section of $S^2(\cO_S\Dsum\cL)$. A local
calculation shows that that section is $\tensor 1\oplus -\varphi(1)$.

We shall have need of the following technical result on Cartier divisors.
%;;%%%Divisibility of Cartier divisors Proposition Cartier divisibility
\begin{proposition}\label{Cartier divisibility}
Let \map{\pi}{X}{S} be a smooth, proper map of schemes with connected fibres and
$D \subset X$ a relative (wrt to $\pi$) effective Cartier divisor and let $n$ be
an integer such that $n\cO_S=\cO_S$.

The functor that to a $T \to S$ associates the set of relative effective Cartier
divisors $D \subset X\times_ST$ such that $nD=C\times_ST$ is representable by a
closed subscheme of $S$. In particular, if $S$ is reduced and for every
geometric fibre of $X \to S$, $C$ is $n$ times an effective Cartier divisor, $C$
itself is $n$ times an effective Cartier divisor.
\begin{proof}
The functor is clearly representable by some $S$-scheme $V \to S$ locally of
finite type and we may, by a standard limit argument, assume that $S$ is
noetherian. What needs to be proven is that $V \to S$ is proper, injective on
geometric points and unramified. For the properness we may use the valuative
criterion so that $S$ is the spectrum of a discrete valuation ring and we assume
given a $D_\eta$ over the generic point with $nD_\eta=nC_\eta$. We then let $D$
be the schematic closure of $D_\eta$ which is a Cartier divisor as $X$ is
regular and is relative as it does not have any horisontal components.

As for injectivity on geometric points we may assume that $S$ is the spectrum of
an algebraically closed field and then the uniqueness of $D$ is clear as $X$ is
regular and thus the group of Cartier divisors is torsion free.

Finally, to prove that $V \to S$ is unramified it is enough to show that it is
formally unramified so we may assume that $S=\Sp R$, where $(R,\maxid_R)$ is a
local Artinian ring, $0\ne\delta \in R$ annihilates $\maxid$ and we assume that
a $D$ exists over the closed subscheme defined by $\delta$. We then are to prove
that there is at most one lifting of $D$ to $X \to S$. Now, a Cartier divisor
$E$ is given by specifying a line bundle $\cL$ and an injective
$\cO$-homomorphism $\cO \to \cL$. As $X \to S$ is flat, the injectivity follows
from injectivity over the special fibre and is hence automatic in our
situation. Furthermore, if $D$ is given by \map{s}{\cO}{\cL}, then $nD$ is given
by \map{s^{\tensor n}}{\cO}{\cL^{\tensor n}}. In our situation we assume a pair
$(\cM,t)$ over $X$ representing $C$ and two pairs $(\cL,s)$ and $(\cL',s')$
whose $n$'th powers are isomorphic to $(\cM,t)$ and whose reductions modulo
$\delta$ are isomorphic. Now the kernel and cokernel of the reduction $\Pic(X)
\to \Pic(X/\delta)$ are $\cO_S$-modules so that multiplication by $n$ is by
assumption bijective on them which shows that $\cL'$ and $\cL$ are isomorphic
and we may assume them to be equal. Hence $s'$ is of the form $s+\delta w$ with
$w$ a section of $\overline{\cL}=\cL/\maxid$ by the flatness of $X \to S$. By
assumption their $n$'th powers are isomorphic so that $(1+\delta \lambda)
s^n=(s+w)^n=s^n+ns^{n-1}\delta w$ for some $\lambda \in R$. This gives
$\overline{s}^{n-1}(\lambda \overline{s}+nw)=0$, where $\overline{s}$ is the
reduction of $s$ modulo $\maxid$. As $\overline s$ is a non-zero divisor this
gives $w =-\lambda/n\overline s$, i.e., $\delta w=-\lambda/n\delta s$ which
gives $s'=(1-\lambda/n\delta)s$ so that the pairs $(\cL,s)$ and $(\cL',s')$ are
isomorphic.

The last statement follows immediately from the previous ones as under its
assumptions $V$ has the same topological space as $S$.
\end{proof}
\end{proposition}
The following result is no doubt well known but I do not know of a reference.
%;;%%%Deformation of curve maps Proposition curve map deformations
\begin{proposition}\label{curve map deformations}
\part Let \map{f}{C}{D} be a separable non-constant map of smooth proper curves
over a field $\k$ and consider the deformation functor that whose values on a
nil-thickening of $\Sp \k$ are isomorphism classes of deformations of $C$
and the map $f$. The map that maps such deformations to similar deformations of
the formal completions of $C$ resp.\ $D$ along the ramification resp.\ branch
locus is an isomorphism.

\part Let $\k$ be a field, $n$ an integer invertible in $\k$ and
\map{f}{\Spf\pow{x}}{\Spf\pow{t}} be the map $t \mapsto x^n$. Then $t \mapsto
x^n+\sum_{0\le i<n-1}a_ix^i$, where the $a_i$ are power series variables, is a
miniversal deformation of $f$.
\begin{proof}
The first part can be proved by noticing that outside of the
ramification/branch loci the map is unramified and hence extends uniquely along
any nil-thickening. This shows that the deformation problem is the same as that
for the localisation along the ramification/branch loci. The comparison between
the deformation for the localisations and the completions is also clear as when
one inverts generators for the ramification/branch loci then the map is
\'etale. This means that the map is specified by choosing a lattice in the ring
of functions in the source over the ring of functions of the target. Specifying
such a lattice is the same in the localisation as in the completion.

Alternatively one can use deformation theory. If $R \to S$ is a small extension
of local Artinian algebras with residue field $\k$, small meaning that the
kernel $K$ is killed by the maximal ideal of $R$, then the liftings of a
deformation over $S$ to one over $R$ is in bijection with
$H^0(C,f^*T_D/T_C)\Tensor_{\k} K$. Indeed, if the deformation of $C$ is kept
fixed then liftings of deformations of $f$, given one, are in bijection with
$H^0(C,f^*T_D)\Tensor K$. Taking into account the possibility of varying also
deformations of $C$ we have to divide out by the action of liftings of
automorphisms of the deformation of $C$, i.e., sections of $T_C\Tensor K$. This
action is given by addition composed with the map $T_C \to f^*T_D$ and hence the
full problem is in bijection with $H^0(C,f^*T_D/T_C)\Tensor_{\k} K$. As we never
used the properness the same is true for the local or complete problem as
$f^*T_D/T_C$ is supported on the ramification locus.

As for the last part, the formula $t \mapsto x^n+\sum_{0\le i<n-1}a_ix^i$ gives
a deformation over $\pow[\k]{a_0,\dots,a_{n-2}}$ and hence a map to the
miniversal deformation. As $\pow[\k]{a_0,\dots,a_{n-2}}$ is (formally) smooth,
to show that this map is an isomorphism it is enough to show that induces an
isomorphism on tangent spaces and for that we can use the description of
deformations over $\k[\delta]$, $\delta^2=0$, just given to show that. Indeed,
the action of the sections of $T_C$ on such deformations is by interpreting a
derivation of $\pow[\k]{x}$ as an automorphism of the scalar extension to
$\k[\delta]$, $\pow[\k]{x}[\delta]$, that is the identity modulo $\delta$ and
then composing the given map $\pow[\k]{t}[\delta] \to \pow[\k]{x}[\delta]$ with
that automorphism. If the vector field is $h(x)d/dx$ and the map has the form
$f(x)+g(x)\delta$ with $f,g,h \in \pow[\k]{x}$ then this composite is
$f(x)+(h(x)f'(x)+g(x))\delta$. This shows that the tangent vector of the map is
given by the residue of $g$ modulo $f'(x)$ which makes it clear that the tangent
map is an isomorphism.
\end{proof}
\end{proposition}
When we make an infinitesimal study of the moduli stack we shall not just deal with
the stack as such but also with the natural stratification of it given by the
ramification exponents of a map between curves. We recall its definition and
first properties given in \cite[App.]{ekedahl01::hurwit}. We begin by noting
that for technical reasons we shall need to assume that we deal with schemes
over $\Sp\Q$ for the rest of this section.
\begin{remark}
Note that this restriction is not just due to the fact that one would need some
slight modifications to get similar results in positive characteristic. In fact
there are some truly new phenomena in positive characteristic. Consider for
instance the case of Proposition \ref{primary}. Condensed it says that for a
finite flat map there is a stratification of the base such that on each stratum
there is a closed subscheme of the total space which is \'etale over the base
and whose defining ideal is nilpotent. A similar result is not possible in
positive characteristic. Consider for instance an inseparable field extension $k
\subset K$ of degree $p$ say. For the corresponding map of schemes $\Sp K \to
\Sp k$ if it had a similar stratification then there could only be one stratum
but $\Sp K$ does not have a closed subscheme which is \'etale over $\Sp k$ and
whose ideal is nilpotent.
\end{remark}
Recall that if \map{f}{Y}{X} is a finite flat map then we define its
\Definition{trace form} to be the symmetric bilinear form $(r,s) \mapsto
\Tr(rs)$. We then define, for each natural number $n$, the closed subscheme of
$Y$ given by the condition that the corank of the map $f_*\cO_Y \to
\Hom_{\cO_X}(f_*\cO_Y,\cO_X)$ induced by the trace form is $\ge n$. We shall
call the stratification thus obtained the \Definition{trace stratification} wrt
$f$. In an open stratum we get the following primary decomposition result.
%;;%%%Primary decomposition Proposition primary
\begin{proposition}\label{primary}
Let \map{f}{Y}{X} be a finite flat map for which $X$ equals a single open trace
stratum. Then the radical of the trace form (i.e., the kernel of the map
$f_*\cO_Y \to \Hom_{\cO_X}(f_*\cO_Y,\cO_X)$) is a subbundle and an ideal. The
closed subscheme $Y'$ defined by it is an \'etale covering of $X$. Furthermore,
there is a unique $X$-retraction $Y \to Y'$ which makes $Y$ a flat $Y'$-scheme.

Finally, if $Y \to Y'$ has rank $n$ then the $n$'th power of the radical of the
trace map is zero.
\begin{proof}
The fact that the radical is a subbundle follows directly from the fact that by
assumption $f_*\cO_Y \to \Hom_{\cO_X}(f_*\cO_Y,\cO_X)$ has constant rank (in the
schematic sense defined by the vanishing and non-vanishing of subdeterminants
when $X$ is not reduced) and then its kernel is a subbundle. That it is an ideal
follows directly from the definition of the trace form. Then $Y' \to X$ is flat
so to prove that it is \'etale it is enough to do it when $X$ is the spectrum of
a field in which case it is well known and easy to see that $Y'$ is
\'etale. Assume that we know that the existence and unicity of a retraction
locally.  Then the unicity forces it to exist globally by (\'etale) descent. By
descent again, the flatness of the retraction needs to be checked only in the
case when $Y'$ is the disjoint union of copies of $X$ in which case it is clear.

For the last statement, replacing $X$ by $Y'$ we may assume that $f$ has rank
$n$. If $\cI$ is the radical and $x \in \cI$ (i.e., is a local section of that
sheaf) then we have that $x^i \in \cI$ for all $i > 0$ and hence
$\Tr(x^i)=0$. As we are in characteristic zero this implies that the
characteristic polynomial of multiplication by $x$ is $t^n$ and by the
Cayley-Hamilton theorem $x^n=0$. Again as we are in characteristic zero we get
by polarisation that any product of $n$ local sections of $\cI$ is zero.

It remains to prove local existence and unicity of a retraction. Locally we may
assume that $Y'$ is a disjoint union of copies of $X$ which makes the existence
of a retraction obvious and the unicity clear.
\end{proof}
\end{proposition}
When the corank of the trace form is constant, the proposition shows that the
function on the points of $Y$ defined by the rank at a point of the radical is a
locally constant function and hence the function which to a point of $X$
associates the multisets of those ranks is locally constant. The stratification
obtained in this way will be referred to as the \Definition{stratification by
multiplicity}. Thus while the trace stratification is a decreasing sequence of
closed subschemes, the stratification by multiplicity is a further decomposition
of the open strata. We shall also in this case denote the subscheme defined by
the radical by $X^{fred}$ and call it the \Definition{fibrewise reduced subscheme}.

Using the primary decomposition we get an extension of the pointwise result that
expresses an effective Cartier divisor on a smooth curve as the sum of points.
%;;%%%Fred of Cartier divisor Proposition Cartier fred
\begin{proposition}\label{Cartier fred}
Let \map{f}{X}{S} be a smooth map of $\Sp\Q$-schemes of relative dimension $1$
and suppose $D \subset X$ is a relative effective Cartier divisor and assume
that the corank of the trace form of $D \to S$ is of locally constant rank. 

\part The fibrewise reduced subscheme $D^{fred}$ of $\Sp\cO_x/\cI_D$ is a
relative effective Cartier divisor. It can be written as the disjoint of
subschemes $D_i$ having the property that at a point $d \in D_i$ the defining
ideal of $D_i$ in $D$ has rank $e_i$ as $\cO_{D_i}$-module. If that is done then
we have that $D=\sum_ie_iD_i$ as divisors.

\part Conversely if $D$ can be written as $D=\sum_ie_iD_i$ with $D_i$ \'etale
disjoint Cartier divisors then the corank of the trace form of $D \to S$ is
locally constant and the union of the $D_i$ is the fibrewise reduced subscheme
of $D$.
\begin{proof}
That $D^{fred}$ is a Cartier divisor is clear and we get from proposition
\ref{primary} that its defining ideal $\cI$ in $D$ is a locally free
$\cO_{D^{fred}}$ and hence its rank is locally constant. This gives the
components $D_i$ and to prove the equality of Cartier divisors we may work
locally around one of the $D_i$, i.e., assume that the rank of $\cI$ as
$\cO_{D^{fred}}$-module is everywhere equal to some $n$. By \(primary) the
$n$'th power of $\cI_{C^{fred}}$ is contained in $\cI_D$. To check that it is an
equality it is enough to check on fibres over closed points of $S$ and then it
is true as they have the same degree at all the points of $D^{fred}$.

Finally, if $D=\sum_ie_iD_i$ and $D'$ is the union of the $D_i$ then
$\cI_{D'}/\cI_D$ is an $S$-flat nilpotent ideal of $\cO_X/\cI_D$ such that the
quotient by it is \'etale. This shows that $\cI_{D'}/\cI_D$ is the radical of
the trace form and thus that $D'$ is the fibrewise reduced subscheme of $D$.
\end{proof}
\end{proposition}
The multiplicities $e_i$ can be considered as locally constant functions on
$D^{fred}$ and we may choose the $D_i$ such that the $e_i$ are all
distinct. Having done that the decomposition is unique and we shall call it the
\Definition{primitive decomposition} of $D$. We shall also consider the locally
constant function on the base $S$ which to a point $s$ associates the multiset
of the multiplicities of the points of $D$ in the fibre over $s$. This will be
called the \Definition{multiplicity multiset} associated to $D$.

If $S$ is a scheme and \map{f}{Y}{X} is a finite $S$-map between smooth
(possibly formal) $S$-schemes of relative dimension $1$ then we get two finite
$S$-schemes, the ramification locus which is a relative Cartier divisor of $Y$
and the branch locus which is a relative Cartier divisor of $X$ (by definition
the branch locus is the norm wrt to $f$ of the ramification locus considered as
a Cartier divisor which is defined as $f$ is finite flat). We shall call the
stratifications by multiplicity on $S$ induced by them the \Definition{ramification
stratification} and \Definition{branch stratification} respectively.
\end{section}
%;;%%Relations with Abel
\begin{section}{The original problem}

We shall begin by formulating in modern terms Abel's question and the answer he
gave to it. The initial setup is that of a square free monic polynomial $R(x)$
of even degree over the complex numbers and the rational differential form
$\omega := \rho dx/\sqrt{R}$, $\rho \in \C(x)$, on the compact Riemann surface
$C$ with field of rational functions $\C(x,\sqrt R):=\C(x)[y]/(y^2-R(x))$. The
general question Abel poses is when this form is the logarithmic differential
$\dlog f:=df/f$ for a non-zero rational function $f$. If $\iota$ is the
hyperelliptic involution of $C$ which takes $x$ to $x$ and $\sqrt R$ to $-\sqrt
R$, then $\iota^*\omega=-\omega$ and as $f$ is determined up to a constant we
get $f\circ\iota\cdot f=\lambda \in \C^*$. Modifying $f$ by multiplying by a
square root of $\lambda^{-1}$ allows us to assume that $f\circ\iota\cdot
f=1$. By a somewhat anachronistic appeal to Hilbert's Theorem 90 we get that $f$
has the form $g/g\circ\iota$ and by clearing denominators we may assume that
$g=P+\sqrt{R}Q$ with $P,Q \in \C(x)$ which is indeed the form that Abel assumes
the solution to have. We note that $g$ is uniquely determined up to a rational
function in $x$. Abel then almost immediately restricts himself to the case
where $\rho$ is a polynomial. This implies that $\omega$ is regular over
$C^\circ:=\Sp\C[x,\sqrt R]$ and hence in particular that $f$ does not have poles
or zeroes in $C^\circ$ or otherwise put, if $\infty_1$ and $\infty_2$ are the
two points of $C$ in the complement of $C^\circ$ in $C$, then
$(f)=m\infty_1-m\infty_2$ for some integer $m$. For a divisor $D$ on $C$ we
denote by $D^\circ$ the part of $D$ that has support on $C^\circ$ and then we
have $0=(f)^\circ=((g)-\iota(g))^\circ$. Now by assumption $(g)^\circ \ge 0$ and
hence for $r \in C^\circ$, $r$ and $\iota r$ appear with the same multiplicity
in $g$. Now, for any $r \in C^\circ$ $r+\iota r-(\infty_1+\infty_2)$ is
the divisor of a rational function in $x$ so that if $r$ is a non-Weierstrass
point (i.e., $r \ne \iota r$) we may modify $g$ by a rational function in $x$ so
that $r$ does not appear in $(g)$. Similarly, if $r$ is a Weierstrass point we
may assume that it appears with at most multiplicity $1$ in $(g)$. In
particular, $(g)^\circ$ has support at the Weierstrass points of $C$.

Even though Abel treats the general case we shall only be interested in the case
when $(g)^\circ=0$. The reason for this is that we shall be mainly interested in
the existence of a $g$ as a condition on the curve $C$ and we have that
$2\omega=\dlog f^2=\dlog g^2/g^2\circ\iota$ and all the Weierstrass points
appear with even multiplicity in $(g^2)$ and they can therefore be removed
completely. Hence at the price of possibly replacing $\rho$ with $2\rho$ we see
that Abel's problem has been reduced to the problem of finding $g \in \C(x,\sqrt
R)$ with $(g)^\circ=0$, i.e., $(g)=n(\infty_1-\infty_2)$ for some integer $n$
and excluding the trivial case of $\rho=0$ we may assume that $n$ is
non-zero. In any case $(g)^\circ=0$ and hence $g$ is a unit in $\C[x,\sqrt
R]$. This implies that $g=P+\sqrt{R}Q$, with $P, Q \in \C[x]$ and furthermore
that the norm $N(g)=P^2-RQ^2$ of $g$ with respect to the finite flat extension
$\C[x,\sqrt R]/\C[x]$ is a unit in $\C[x]$, i.e., a non-zero constant. After
changing $g$ by a constant we may assume that $P^2-RQ^2=1$. Abel then notes that
this is analogous to Pell's equation and proceeds to use continued fractions in
analogy with the case of Pell's equation. There is a difference however in that
the number theoretic case gives a method for solving Pell's equation while the
geometric case gives a criterion for the existence of a solution (as well as a
method for constructing it when it does exist).
\begin{remark}
Abel's approach gives a condition on $R$ for a solution to the problem with $n$
arbitrary to exist. This is not appropriate for our purposes as when $n$ varies
we get a countable union of closed subvarities in the space of $R$'s which is
unnatural from a geometric point of view.
\end{remark}
\end{section}

%;;%%Continued fraction

%;;%%Abel curves sec:Abel curves
\begin{section}{Abel curves}\label{sec:Abel curves}

We shall now give the formal definition of an Abel curve. To simplify the
presentation (and make it closer to Abel's original results) we shall \emph{from
now on assume that $2$ is invertible in all our schemes}. If we want to consider
the moduli space of Abel curves we want to make as few choices as possible as
any choice leads to a larger space which is the reason for the somewhat lengthy
definition. To avoid ambiguities in the case of genus $1$ to us a
\Definition{hyperelliptic curve} will be a smooth proper curve $C$ together with
a \emph{choice} of an involution $\iota$ such that the quotient $C/\iota$ is of
genus zero.\footnote{With this definition we can have hyperelliptic curves of
genus zero which for our purposes is quite acceptable though rather trivial.}
%;;%%%Abel curves Definition Abel curves
\begin{definition}\label{Abel curves}
A \Definition{(smooth) Abel curve of genus $g$ and order $n$} over a scheme $S$
consists of
\begin{itemize}
\item a smooth and proper $S$-curve \map{\pi}{C}{S},

\item an $S$-involution $\iota$ of $C$ making each fibre a hyperelliptic curve
of genus $g$, 

\item two disjoint sections $\infty_1$ and $\infty_2$ of $\pi$ such that
$\infty_2=\iota\infty_1$, 

\item a line bundle $\cL$ on $S$ together with a trivialisation $\varphi\co\cO_S
\riso \cL^2$, and
%%%%%%%%%%%$

\item a finite flat $S$-map \map{f}{C}{P}, where $(P,\sigma)$ is the involutive
bundle associated to $\cL$ and the trivialisation $\varphi$, of degree $n$ such
that the sections $\infty_1$ and $\infty_2$ map to the sections $\mathbf{0}$ and
$\infty$ of the involutive bundle $\P$ and $\sigma\circ f=f\circ\iota$.
\end{itemize} 
A \Definition{split Abel curve} is an Abel curve together with a splitting of
the involutive bundle.

An isomorphism between Abel curves consists of isomorphisms between the $C$ and $P$
parts of the curves transporting all the structures of the first curve to those
of the second.

Associating to each $S$ the groupoid of Abel curves and isomorphisms
between them gives a stack (in say the flat topology) that we shall denote
$\Ab{g}{n}$ and similarly we get the stack of split Abel curves $\sAb{g}{n}$.
\end{definition}
When the base is an algebraically closed field we get exactly the description
that came out of Abel's problem. Note that in that case it follows from the
equation $P^2-RQ^2=1$ that $2\deg P \ge \deg R$, i.e., $2n \ge 2g+2$ which means
$n \ge g+1$.

It is not immediately clear that this is the right definition for families as
one could worry that we have made an unnecessary choice in choosing two sections
$\infty_1$ and $\infty_2$ instead of a divisor of degree $2$ that only after a
base change splits up into two disjoint sections. The following definition
expresses that concern.
%;;%%%Abel curves Definition twisted Abel
\begin{definition}\label{twisted Abel}
A \Definition{twisted (smooth) Abel curve of order $n$} over a scheme $S$
consists of
\begin{itemize}
\item a smooth and proper $S$-curve \map{\pi}{C}{S},

\item an $S$-involution $\iota$ of $C$ making each fibre a hyperelliptic curve, 

\item a $\iota$-invariant relative effective divisor $D$ of degree $2$ of $C$
which is \'etale over $S$ and on which $\iota$ acts freely,

\item a smooth and proper $S$-curve \map{\rho}{P}{D} all of whose fibres have
genus zero,

\item an $S$-involution $\sigma$ of $P$,

\item a $\sigma$-invariant relative effective divisor $D'$ of $P$ which is
\'etale over $S$ and on which $\sigma$ acts freely,

\item a finite flat $S$-map \map{f}{C}{P} of degree $n$ such that the inverse
image (as effective divisors or equivalently as subschemes) of $D'$ is $nD$ and
for which $f\circ\iota=\sigma\circ f$.
\end{itemize}
Isomorphisms between Abel curves consist of an isomorphism $g$ between the
$C$-parts preserving the $\iota$'s, $\infty_1$'s, and $\infty_2$'s and an
automorphism $h$ of $\P^1\times S$ such that $f\circ g=h\circ f$.

Associating to each $S$ the groupoid of twisted Abel curves and isomorphisms
between them gives a stack (in say the flat topology) that we shall denote
$\tAb{g}{n}$.
\end{definition}
The relation between these definitions is expressed in the following result.
%;;%%%Twisted stack Proposition twisted reduction(1)
\begin{proposition}
\part \label{twisted reduction} The stack of twisted Abel curves of genus $g$
and order $n$, $\tAb{g}{n}$, is equivalent to $B\Sigma_2\times\Ab{g}{n}$, where
$B\Sigma_2$ is the stack of $\Sigma_2$-torsors, i.e., the stack of \'etale
double covers. The projection on the first factor associates to a twisted Abel
curve, using the notation of Definition \ref{twisted Abel}, the \'etale double
cover $D \to S$.

\part The forgetful map $\sAb{g}{n} \to \Ab{g}{n}$ is an \'etale double
cover. 
%% For an Abel curve $(C \to S,\iota,\infty_1,\infty_2,f,P)$, $\cL$ and the
%% trivialisation $\cO_S \riso \cL^2$ is recovered by its associated double 
%% cover being the involutive locus.
\begin{proof}
Using the notation of Definition \ref{twisted Abel} we get from a twisted Abel
curve over $S$ an \'etale double cover $D \to S$ which gives a map from the
stack of twisted Abel curves to $B\Sigma_2$. On the other hand, $(\iota,\sigma)$
gives an involution of the Abel curve and we may use it and the double cover $D
\to S$ to twist the Abel curve, in particular the twist, $\tilde C$, of $C$ is
obtained by taking the quotient of $D\times_SC$ by the action of
$(\iota,\iota)$. The section given by the graph of the inclusion of $D$ in $C$
is invariant under this map and hence descends to a section of $\tilde C \to S$
and the same is true of the group of the map $D \to C$ composed with $\iota$. In
other words, the divisor $D$ twists to give a divisor that is the disjoint union
of two sections. Now, the map \map{f}{C}{P} maps $D$ isomorphically to $D'$ so
that also $D'$ is the disjoint union of two sections. The existence of these two
disjoint sections makes $P \to S$ isomorphic to $\P(\cL\Dsum\cM)$ for some line
bundles $\cL$ and $\cM$ on $S$, where the two sections correspond to the two
summands. Now, $\sigma$ permutes the two section which forces $\cL$ and $\cM$ to
be isomorphic so that $P \to S$ is isomorphic to $\P^1\times S \to S$ with the
two sections given by $0\times S$ and $\infty\times S$. As the inverse images of
$0\times S$ and $\infty\times S$ are $n$ times the two sections of $D$ we get an
Abel curve and consequently a map $\tAb{g}{n} \to \Ab{g}{n}$ and combining the
two constructed maps we get a map $\tAb{g}{n} \to
B\Sigma_2\times\Ab{g}{n}$. Conversely, given an Abel curve over $S$ we can
consider the map $f$ as an isomorphism $\cO_C(n\infty_1) \riso
\cO_C(n\infty_2)$. Letting $\iota$ act on that isomorphism gives another
isomorphism $\cO_C(n\infty_2) \riso \cO_C(n\infty_1)$. Their composites are
then multiplication by an invertible function $\lambda$ on $S$. That means that
if we define $\sigma$ on $\P^1\times S$ by $(x:y) \mapsto (\lambda y:x)$ then
$f\circ\iota=\sigma\circ f$ so that we have a twisted Abel curve over $S$. Now,
$(\iota,\sigma)$ is an involution of that object and so that if we have an
\'etale double cover $\tilde D \to S$ we can use it to twist our twisted Abel
curve and we obtain thus a map $B\Sigma_2\times\Ab{g}{n} \to \tAb{g}{n}$ which
is clearly an inverse to the map just constructed.

As for the second part it is clear.
\end{proof}
\end{proposition}
The proposition shows that it is no real loss in generality to restrict
ourselves to Abel curves which we shall do that from now on with the exception
of the following result which confirms the representability of the two stacks.
%;;%%%Representability
\begin{proposition}
The stacks $\Ab{g}{n}$, $\sAb{g}{n}$, and $\tAb{g}{n}$ are Deligne-Mumford stacks of finite
type over $\Sp\Z[1/2]$.
\begin{proof}
This is quite standard as soon as we have verified that the automorphism group
scheme of an (twisted) Abel curve over an algebraically closed field is finite
\'etale. For $g \ge 2$ this is clear as it is true for all curves of genus
$g$. For $g=1$ we have to use the fact that the hyperelliptic involution is part
of the structure so that an automorphism has to commute with it. For a
hyperelliptic involution $\iota$ we may choose a fix point as origin and in the
thus obtained group structure on the curve, the involution is multiplication by
$-1$ and then it is clear that the automorphism group scheme centralizer of the
curve is finite \'etale. Finally, for genus $0$ we have to look at the
automorphism group scheme of automorphisms of $\P^1$ fixing two points and
commuting with an involution that permutes the two points. It is clear that the
points and the involution is conjugate to $0$, $\infty$, and $x \mapsto 1/x$ and
then the automorphism group scheme that fixes these is clearly finite \'etale.
\end{proof}
\end{proposition}
Our definition of an Abel curve is chosen to be closely modeled on Abel's
original condition. On the other hand -- at least punctually -- the relevant
condition is that the divisor class $\infty_1-\infty_2$ is killed by $n$ as
then there is a map to $\P^1$ with whose zero and pole divisor is
$n(\infty_1-\infty_2)$. This turns out to be true for families.
%;;%%%Jacobian section Proposition Jacobian
\begin{proposition}\label{Jacobian}
Let $\Hyp{g}{2}$ be the stack of hyperelliptic curves with two distinct points
$(C,\iota,a,b)$ of genus $g$ and let \map{s}{\Hyp{g}{2}}{J_g} be the section of
the Jacobian of the universal curve given by $a-b$. Let $\cH$ be the closed
substack of $\Hyp{g}{2}$ defined by the conditions $\iota a=b$ and $ns=0$. Let
$\rho$ be the involution of $\Ab{g}{n}$ which takes an object $(C \to
S,\iota,\infty_1,\infty_2,f)$ to $(C \to S,\iota,\infty_1,\infty_2,-f)$. Then
the map given by
\begin{displaymath}
\function{\Ab{g}{n}}{\cH}{(C,\infty_1,\infty_2,f)}{(C,\infty_1,\infty_2)}
\end{displaymath}
is an isomorphism of stacks.
\begin{proof}
As has been noted above, $f$ may be thought of as an isomorphism
$\phi\co\cO_C(n\infty_1) \riso \cO_C(n\infty_2)$ and then $\rho$ takes it to
$-\phi$. On the other hand, an $S$-object of $\cH$ has the property that
$\cO(n\infty_1-n\infty_2)$ is a pullback of a (unique) line bundle $\cL$ on
$S$. Now, applying $\iota$ to $\cO(n\infty_1-n\infty_2)$ gives its inverse
which translates into an isomorphism $\cL \riso \cL^{-1}$, i.e., a
trivialisation of $\cL^2$. This gives an object of $\Ab{g}{n}$ over $S$.
\end{proof}
\end{proposition}
As the zero-section in an abelian scheme is a complete intersection subscheme we
get one immediate consequence.
%;;%%%Lower dimension bound.
\begin{corollary}
$\Ab{g}{n} \to \Sp\Z[1/2]$ is of relative dimension at least $g$ at each of its
points and at a point where the relative dimension is $g$ it is a local complete
intersection.
\begin{proof}
The substack $\cH$ of $\Hyp{g}{2}$ fulfilling $\iota a=b$ is an open substack of
the stack $\Hyp{g}{1}$ of hyperelliptic curves with one chosen point, namely the
complement of the locus of fixed points of the hyperelliptic involution, where
the isomorphism maps $(C,\iota,a)$ to $(C,\iota,a,\iota a)$. Hence that substack
is smooth of relative dimension $2g-1+1=2g$. Now, by the proposition $\Ab{g}{n}$
is the inverse image in $\cH$ of the zero section of $J_g \to \Hyp{g}{2}$ under
the map $ns$ and $J_g \to \Hyp{g}{2}$ being smooth, the zero section is a local
complete intersection map of codimension $g$.
\end{proof}
\end{corollary}
\begin{remark}
In characteristic $0$ we shall show that the codimension is in fact $g$ and that
$\Ab{g}{n}$ is in fact smooth.
\end{remark}
Fix $n$ and $g$ with $n \ge g+1$ and consider $\A:=\A^{2n+g+3}_{\Z[1/2]}$ that
we shall regard as the parameter space for triples $(P,Q,R)$ of polynomials of
degrees $n$, $n-g-1$, and $2g+2$ respectively with $R$ monic. We let
\Definition{$\vAb{g}{n}$} be the subscheme of triples that fulfill $P^2-RQ^2=1$
and for which $R$ is square free (i.e., its discriminant is invertible) and $P$
and $Q$ have invertible top coefficients. We let $\uAb{g}{n}$ be the subscheme
of $\vAb{g}{n}$ defined by the condition that $R$ is normalised and $P$ and $Q$
are monic. Over $\vAb{g}{n}$ we have an Abel curve given by
\begin{displaymath}
 C:=\Proj\cO_{\uAb{g}{n}}[s,t,y]/(y^2-t^{2g+2}R(s/t)),
\end{displaymath}
where $\deg s=\deg t=1$ and $\deg y=g+1$, $\iota$ is given by $(s\co t\co y)
\mapsto (s\co t\co -y)$, $\infty_1$ and $\infty_2$ are given by $(0\co1\co1)$
resp.\ $(0\co1\co-1)$, and $f$ is given by $(s\co t\co y)\mapsto
(t^n(P(s/t)+yQ(s/t))\co t^n)$. This therefore gives a map $\vAb{g}{n} \to
\Ab{g}{n}$. We shall call any Abel curve that is a pullback of this family by a
map to $\vAb{g}{n}$ a \Definition{Pell family} and if it is given as a pullback
by a map to the closed subscheme $\uAb{g}{n}$ we shall call it a
\Definition{normalised Pell family}.
%;;%%%Pell versality Theorem Pell versality
\begin{theorem}\label{Pell versality}
\part $\vAb{g}{n} \to \Ab{g}{n}$ factors through the map $\sAb{g}{n} \to
\Ab{g}{n}$.

\part[v part] Over $\Sp\Q$ the map $\vAb{g}{n} \to \sAb{g}{n}$ is a torsor under
the subgroup of $\mul\times\Aff$, where $\Aff$ is the group of affine
transformations of the affine line, of pairs $(\lambda,z \mapsto az+b)$ for
which $\lambda^2=a^{2g+2}$.

\part[u part] Over $\Sp\Q$ the map $\uAb{g}{n} \to \sAb{g}{n}$ is a torsor under
the subgroup of $\mul\times\Aff$ of pairs $(\lambda,z \mapsto az+b)$ for which
$\lambda=a^{g+1}$, $a^n=1$ and $b=0$, a group isomorphic to the group $\mu_n$ of
$n$'th roots of unity.

\part In particular the map $\vAb{g}{n} \to \Ab{g}{n}$ is a
chart. i.e., smooth and surjective, and $\uAb{g}{n} \to \Ab{g}{n}$ is even an
\'etale chart.
\begin{proof}
To prove the first part we note that for a Pell family the involution on
$\P^1\times S$ compatible with $f$ and $\iota$ is $x \mapsto 1/x$ whose fixed
point scheme is $\pm1$ and by ordering it as $\{1,-1\}$ we get a family in
$\sAb{g}{n}$. 

Assume now that $(C \to S,\iota,\infty_1,\infty_2,f)$ is a family in
$\sAb{g}{n}$. By assumption, using the notation of definition \ref{Abel curves},
$P$ is isomorphic to $\P^1\times S$ in a way such that $\mathbf{0}$ on $P$ is
$0\times S$ and $\infty$ is $\infty\times S$ and the involution $\sigma$ is $x
\mapsto 1/x$. Consider now the quotient $D$ of $C$ by $\iota$. As $2$ is
invertible taking the quotient by $\iota$ commutes with base change so that in
particular \map{\pi}{D}{S} is a smooth proper map with genus $0$
fibres. Furthermore, either of the sections $\infty_1$ or $\infty_2$ give a
section $\infty$ of $\pi$. Now, again as $2$ is invertible, the double cover $C
\to D$ is given by a line bundle $\cM$ on $D$ and a section of $\cM^2$. As $\cM$
has degree $g+1$ on each fibre $\cM(-(g+1)\infty)$ is the pullback from $S$ of
a line bundle $\cL$.

We shall now show that giving an isomorphism of $D$ with $\P^1\times
S$ taking $\infty$ to $\infty\times S$ and trivialising $\cL$ is the
same thing as giving a Pell family over $S$ and an isomorphism with it
and our split Abel curve. This will prove the second part and the
third follows as the group of affine transformations is smooth.

In one direction it is clear as a Pell family gives by construction a
trivialisation of $D$ as well as $\cL$. 

For the converse we shall need to use (cf., \(Abel smoothness) which assumes
that we are over $\Q$) that $\uAb{g}{n}$ is smooth so that we may assume that
$S$ is smooth.\footnote{In fact we only use that it is reduced.}

Assume now that an isomorphism $D \riso \P^1\times S$ and a trivialisation of
$\cL$ has been given. This means that $\cM$ is isomorphic to $\cO(g+1)$ so that
the section of $\cM^2$ is a homogeneous form $R(s,t)$ of degree $2g+2$ with
coefficients in $\Gamma(S,\cO_S)$. The existence of the sections $\infty_1$ and
$\infty_2$ show that $R(1,0)$ is a non-zero square and hence after scaling $R$
we can assume that $R(s,1)$ is monic. Now, as the Abel curve is split we may
regard $f$ as an isomorphism \map{f}{\cO_C(\infty_1)}{\cO_C(\infty_2)} and then,
again by the fact that the curve is split, $f\circ\iota^*(f)$ is scalar
multiplication by a square and hence by scaling $f$ we may assume that
$f\circ\iota^*(f)=1$. On $C^\circ:=C\setminus\{\infty_1\}\cup\{\infty_2\}$ $f$
maps into $\mul_S$ so that $f$ is a unit in $\Gamma(C^\circ,\cO)$. This ring is
equal to $\Gamma(S,\cO_S)[s,y]/(y^2-R(s,1))$ so that $f$ has the form
$P(s)+yQ(s)$ and the condition $f\circ\iota^*(f)=1$ translates into
$P^2-RQ^2=1$. Now, if the base is a field it is easy to see that the degree of
$P$ is equal to the degree $n$ and hence, as $S$ is reduced $P$ is of degree $n$
and its top coefficient is a unit. The equation $P^2-RQ^2=1$ and the fact that
$R$ is monic shows that $Q$ has degree $n-g-1$ with invertible top coefficient,
i.e., we have a map to $\vAb{g}{n}$. The possible changes in choices is given by
a scaling factor $\lambda$, which is a unit in $\cO_S$, in the choice of
trivialisation of $\cL$ and an affine transformation $s \mapsto as+b$ where $s
\in \Gamma(S,\cO_S^\times)$ and $b \in \Gamma(S,\cO_S)$. This change takes $y$
to $\lambda y$ and then $(P(s),Q(s),R(s))$ to $(P(as+b),\lambda^{-2}
R(as+b),\lambda Q(as+b))$ so that if we want to keep $R$ monic we need
$\lambda^2=a^{2g+2}$ which shows \DHrefpart{v part}.

Turning to \DHrefpart{u part} we may after an \'etale extension which extracts
an $n$'th root of the top coefficient compose with a change of trivialisation
and affine transformation such that $P$ is monic. As $R$ is also monic this
forces the top coefficient of $Q$ to be $\pm1$ and if $-1$ we may change the
trivialisation by $-1$ to get that $Q$ is also monic. We may then by an
appropriate affine transformation of the form $s \mapsto a+s$ assume that $R$ is
normalised, i.e., we have obtained an $S$-point of $\uAb{g}{n}$. The ambiguities
in our choices are then reduced to a pair $(\lambda,s \mapsto as)$ with $a^n=1$,
$\lambda^{2}=a^{2g+2}$ and $1=\lambda a^{n-g-1}$ conditions which are equivalent
to $a^n=1$ and $\lambda=a^{g+1}$.

The last statement is now clear.
\end{proof}
\end{theorem}
\begin{remark}
Despite the very explicit form of these charts it seems difficult to use them. I
have for instance not been able to show the smoothness of the moduli space using
the Pell equation directly (in the generic case when $R$ and $Q$ have no common
zeros it can be done).
\end{remark}
We may use this result to show that $\sAb{g}{n}$ is almost a scheme by computing
the fixed point sets for the action of subgroups of $\mu_n$ on $\uAb{g}{n}$. For
this we introduce $\wAb{g}{n}$ as the closed subscheme of $\vAb{g}{n}$
consisting of tuples $(P,R,Q)$ if $\vAb{g}{n}$ for which $P$ and $Q$ are monic.
%;;%%%Non-scheme points Proposition fixed points
\begin{proposition}\label{fixed points}
Let $m > 1$ be an integer that divides $n$ so that $\mu_m \subseteq
\mu_n$. Then the fixed point locus for $\mu_m$ acting on $\uAb{g}{n}$ is empty
unless $2g+2 \equiv 0,1 \bmod m$.

\part If $m|g+1$ then the fixed point scheme is of the form
$(p(s^m),r(s^m),q(s^m))$, where $(p(t),r(t),q(t))$ is the universal family of
$\wAb{(g+1)/m-1}{n/m}$.

\part If $m|2g+2$ but $m\not|g+1$ then the fixed point scheme is of the form
$(p(s^m),r(s^m),s^{m/2}q(s^m))$, where $(p(t),r(t)t,q(t))$ is the universal
family $(P,R,Q)$ of $\wAb{(g+1)/m-1/2}{n/m}$ restricted to the closed subscheme
given by $R(0)=0$.

\part if $m|2g+1$ then the fixed point scheme is of the form
$(p(s^m),sr(s^m),s^{(m-1)/2}q(s^m))$, where $(p(t),r(t)t,q(t))$ is the universal
family $(P,R,Q)$ of $\wAb{(2g+1)/(2m)-1/2}{n/m}$ restricted to the closed
subscheme given by $R(0)=0$.
\begin{proof}
If a tuple $(P(s),R(s),Q(s)$ is a point of $\uAb{g}{n}$ and $\zeta$ an $m$'th
root of unity then $\zeta$ takes the tuple to $(P(\zeta s),\zeta^{-2g-2}R(\zeta
s),\zeta^{g+1}Q(\zeta s)$. Hence that the tuple is fixed under $\mu_m$ is
equivalent to $P$, $R$, resp.\ $Q$ being homogeneous of degrees $0$, $2g+2$,
resp.\ $-g-1$, where the grading takes values in $\Z/m\Z$ and $s$ has degree
$1$. This means that for a tuple that is a fixed point, $R(s)$ is of the form
$r(s^m)$. Furthermore, if $k$ is the residue modulo $m$ of $2g+2$ then $s^k$
will be the lowest order non-zero monomial of $R$ and as $R$ does not have any
multiple roots this implies that $k$ is $0$ or $1$. Assume that $m|g+1$.  Then
$R(s)$ has the form $r(s^m)$ and $Q(s)$ has the form $q(s^m)$. Clearly, $p$,
$q$, and $r$ are all monic and as $p^2(s^m)-r(s^m)q^2(s^m)=1$ we get
$p^2(t)-r(t)q^2(t)=1$ so that $(p,r,q)$ gives a family in $\wAb{(g+1)/m-1}{n/m}$
and conversely such a family gives a fixed point $(p(s^m),r(s^m),q(s^m))$ (note
that as $m > 1$ $p(s^m)$ is automatically normalised and that $r(t)$ is
multiplicity free precisely when $r(s^m)$ is). Assume that $m|2g+2$ but
$m\not|g+1$. Then we still have $P(s)=p(s^m)$ and $R(s)=r(s^m)$ but
$Q(s)=s^{m/2}q(s)$ and $P^2(s)-R(s)Q^2(s)=1$ gives $p^2(t)-r(t)tq^2(t)=1$ so
that $(p(t),r(t)t,q(t))$ gives a family in $\wAb{(g+1)/m-1/2}{n/m}$ for which
the $R$-component is $0$ at $0$. Finally if $2g+2 \equiv 1 \bmod m$ we get
$P(s)=p(s^m)$, $R(s)=r(s^m)s$, and $Q(s)=s^{(m-1)/2}q(s)$ which gives
$p^2(t)-r(t)tq^2(t)=1$.
\end{proof}
\end{proposition}
\begin{remark}
By the arguments of the proof of Theorem \ref{Pell versality} (and assuming we
are in characteristic zero) $\wAb{g}{n}$ is isomorphic to $\add\times\uAb{g}{n}$
through affine translations $s \mapsto s+a$ in the polynomial variable. The
subscheme defined by $R(0)=0$ is by the same argument isomorphic to the finite
\'etale cover of $\uAb{g}{n}$ whose $S$-object are $(P,Q,R)$, an $S$-object of
$\uAb{g}{n}$, together with a choice of zero of $R$.
\end{remark}
\end{section}
%;;%%%Genus 0 and 1
\begin{section}{Low genera}

It should come as no surprise that the cases of Abel curves of genus $0$ and $1$
are special and we start by treating them.
%;;%%%%Genus 0 and 1
\begin{proposition}
\part $\Ab{0}{n}$ is isomorphic to $B\Sigma_2$ with universal family having
$\P^1$ as curve with hyperelliptic involution $x \mapsto 1/x$, function
\map{f}{\P^1}{\P^1} given by $x \mapsto x^n$ and involution $\sigma(x)=1/x$. The
mapping to $B\Sigma_2$ giving the isomorphism is given by associating to an Abel
curve the fixed point locus of its hyperelliptic involution.

\part Let $\cA_1 \to \cM_1$ be the universal elliptic curve. Let $\cU$ be the
open substack of the fibre square of $\cA_1 \to \cM_1$ which is the complement
of the diagonal and let \map{\varphi}{\cU}{\cA_1} be the map $(x,y) \mapsto
x-y$. Then $\Ab{1}{n}$ is isomorphic to the inverse image of the kernel of
multiplication by $n$ by $\varphi$
\begin{proof}
Starting with the genus zero case suppose we have a family of Abel curves of
genus zero and degree $n$ $(C \to S,f,\infty_1,\infty_2,\cL,\varphi)$. Then
$\cO(\infty_1-\infty_2)$ is the pullback of a (unique) line bundle $\cM$ on $S$
and the involution $\iota$ induces an isomorphism $\cM \riso \cM^{-1}$ (which
identifies $C$ and $\iota$ with the involutive bundle and involution associated
to the obtained trivialisation of $\cM^2$). Now, f corresponds to an isomorphism
$\cO(n\infty_1) \riso \cO(n\infty_2)\Tensor\cL$, i.e., an isomorphism
$\cM^{\tensor n} \riso \cL$ and the fact that $f\circ\iota = \sigma\circ f$,
where $\sigma$ is involutive involution, implies that $\varphi$ equals the
$n$'th power of the given trivialisation $\cM \riso \cM^{-1}$. This shows the
whole Abel curve is determined by the involutive line bundle $\cM$.

As for the genus $1$ case we start by identifying the closed substack $\H$ of
$\Hyp{1}{2}$ of triples $(\iota,a,b)$ with $\iota a=b$. In fact for any two
disjoint sections $a$ and $b$ of a family of genus $1$ curves there is a unique
hyperelliptic involution that takes $a$ to $b$, namely $x \mapsto -x+a+b$. This
implies is isomorphic to $\cU$ and the rest follows from \(Jacobian).
\end{proof}
\end{proposition}
\end{section}
%;;%%Hurwitz type description
\begin{section}{Hurwitz type description}

If \map{f}{C}{\P^1} is a split Abel curve with hyperelliptic involution $\iota$
then $f\circ\iota=f^{-1}$. The map \map{\tau}{\P^1}{\P^1} given by
$\tau(x)=1/2(x+x^{-1})$ is a quotient map for the action of the involution $x
\mapsto x^{-1}$. We therefore get a commutative diagram
\begin{displaymath}
\CD
C@>>>D:=C/\iota\cr
@V{f}VV @V{g}VV\cr
\P^1@>>{\tau}>\P^1
\endCD
\end{displaymath}
and we see that we may recover $C$ from the map $g$ by taking the normalisation
of its pullback along $\tau$. The map $f$ is then also determined. This gives
the possibility of describing Abel curves in terms of maps of the form $g$. This
is precisely what we are going to do in this section.
%;;%%%Abel map Definition def:Abel map
\begin{definition}\label{def:Abel map}
An \Definition{Abel map} of genus $g$ and degree $n$ over a scheme $S$ consists
of
\begin{itemize}
\item a smooth proper map \map{\pi}{P}{S}, the fibres of which are genus $0$
curves,

\item a section $\infty$ of $\pi$ and an effective Cartier divisor $C$ of $P$
that is \'etale over $S$,

\item an involutive line bundle $(\cL,\varphi)$ over $S$ with \map{\pi}{Q}{S}
the associated projective bundle,

\item an $S$-morphism \map{g}{P}{Q} fibrewise of degree $n$ such
that, $g^*\infty=n\infty$ as Cartier divisors, and

\item a relative effective Cartier divisor $D \subset P$ such that
$g^*\cF=C+2D$, where $\cF$ is the fixed point scheme of $\varphi$ which is an
effective Cartier divisor.
\end{itemize}
A \Definition{split Abel map} is an Abel map together with a splitting of the
involutive bundle.
\begin{remark}
Note that $g$ is flat so that $g^*$ of Cartier divisors is well defined.
\end{remark}
\end{definition}
Given an Abel family $(\map{f}{X}{P},(\cL,\varphi),\iota,\sigma)$ over a scheme
$S$, where $P$ is involutive bundle associated to the involutive line bundle
$(\cL,\varphi)$, we may consider that induced map \map{g}{C/\iota}{P/\sigma}. As
$2$ is invertible taking the quotient by an involution commutes with base change
so that $X/\iota \to S$ is a smooth genus $0$ fibration whereas $P/\sigma$ is
the involutive quotient and hence is isomorphic to $\P(\cO\Dsum\cL)$. Consider
now the induced map $X \to P\times_{P/\sigma}X/\iota$. The composite with it and
the projection $P\times_{P/\sigma}C/\iota \to X/\iota$ is the quotient map and
both $X \to X/\iota$ and $P\times_{P/\sigma}X/\iota \to X/\iota$ are double
covers. As such they are specified by line bundles $\cM$ and $\cN$ and sections
$s$ and $t$ of $\cM^{-2}$ resp.\ $\cN^{-2}$. The map $C \to
P\times_{P/\sigma}C/\iota$ corresponds to a map $\cN \to \cM$ compatible with
the sections of $\cM^{-2}$ and $\cN^{-2}$. The map $\cN \to \cM$ defines a
relative Cartier divisor $D$ as it defines a Cartier divisor on each fibre (over
$S$). Let $C$ be the divisor of $s$ and note that the divisor of $t$ is the
pullback by $g$ of the divisor of the involutive quotient map, i.e., $\cF$
where $\cF$ is the involutive locus. The compatibility between the coverings
then gives that $g^*\cF=C+2D$ and as $g^*\infty=n\infty$ as
$f^*\infty=n\infty_1$ we have an Abel map. Finally, as $X$ is smooth $C$ is
\'etale over $S$. This construction can be reversed.
%;;%%%%Abel maps <-> Abel curves
\begin{proposition}
The stack of Abel curves is isomorphic to the stack of Abel maps.
\begin{proof}
We have just defined a map in one direction. Conversely, assume given an Abel
curve and using the notations of definition \ref{def:Abel map} we recall that
the involutive quotient map is given by $\cO_{\P(\cO\Dsum\cL)}(-1)$ and the
section $\tensor 1\oplus \varphi(1)$ of $\cO_{\P(\cO\Dsum\cL)}(2)$ whose Cartier
divisor is the involutive locus $\cF$. The pullback of it by $g$
is then given by $\cN:=g^*\cO_{\P(\cO\Dsum\cL)}(-1)$ and the Cartier divisor
$g^*\cF$. If we put $\cM := \cN(D)$ then by assumption the section of
$\cN^{-2}$ comes from one of $\cM^{-2}$ and hence gives a double covering $X \to
P$ that maps to the $g$-pullback of the involutive double cover and $X$ is
smooth as $C$ is \'etale and $2$ is invertible. This gives an
inverse map.
\end{proof}
\end{proposition}
In the future we shall pass freely back and forth between Abel maps and Abel
curves.
\end{section}
%;;%%Infinitesimal calculations Section sec:infinitesimal
\begin{section}{Infinitesimal calculations}\label{sec:infinitesimal}

In this section we shall study the deformation theory of Abel maps (and hence of
Abel curves). To avoid problems with wild ramification (and worse still,
inseparable maps) we shall from now on assume that all our schemes and stacks
are over $\Sp\Q$.

If $(p_1(x),p_2(x),\dots,p_n(x))$ is a sequence of monic polynomials over $S$
(i.e., with coefficients in $\Gamma(X,\cO_X)$) then we may put $X$ and $Y$ equal
to $\Sp(\Dsum_i\pow[\cO_S]{x})$ and let $f$ be given by $x \mapsto p_i$ on the
$i$'th component and we shall refer to the ramification and branch
stratifications associated to $f$ as the ramification resp.\ branch
stratifications of the sequence $(p_1(x),p_2(x),\dots,p_n(x))$.
%;;%%%Space of polynomials Definition def:Pol
\begin{definition}\label{def:Pol}
Let $T=(S,S_1,S_2)$ be a sequence of disjoint finite sets and $r$ a function
from $S'$, the disjoint union of the components of $T$ to the positive integers,
for $i=1,2$ let $S^e_i$ and $S^o_i$ be the subsets of $S_i$ where $r$ takes even
resp.\ odd values and set
\begin{displaymath}
n := \sum_{s \in S} (r(s)-1) +  \sum_{s \in S^e_1\disjunion S^e_2} (r(s)/2-1) +
\sum_{s \in S^o_1\disjunion S^o_2} (r(s)-1)/2.
\end{displaymath}
We define $\Pol{T}$ to be the affine space $\A^n$ seen as the parameter space of
tuples $(p_s)_{s \in S'}$ where $p_s$ is a normalised polynomial of degree
$r(s)$ if $s \in S$, a normalised polynomial of degree $r(s)/2$ if $s \in
S^e_1\disjunion S^e_2$ and a monic polynomial of degree $(r(s)-1)/2$ if $s \in
S^o_1\disjunion S^o_2$. Despite this interpretation we shall continue to refer
to the origin as the origin.

To a point $(p_s)$ of $\Pol{T}$ we associate the tuple $(q_s)_{s \in S'}$, where
$q_s=p_s$ if $s \in S$, $q_s=p^2_s$ if $s \in S^e_1\disjunion S^e_2$ and
$q_s(x)=(x-a_s)p^2_s(x)$ with $a_s$ being twice the next to highest coefficient
of $p_s$. (Thus $q_s$ is always a normalised polynomial of degree $r(s)$.) The
ramification and branch stratifications of the sequence $(q_s)$ associated to
the tautological sequence will be referred to as simply the ramification resp.\
branch stratification of $\Pol{T}$.
\end{definition}
We have the following characterisation of the points of a stratum.
%;;%%%Placing of point in stratum Proposition index and corank
\begin{proposition}\label{index and corank}
Let \map{f}{Y}{X} be a finite map of (possibly formal) smooth $1$-dimensional
schemes over a field $\k$. Then the corank of the trace map of a closed point
$s$ of the ramification locus of $f$ is equal to the ramification index at $s$
minus $1$ and the corank of the trace map of a closed point $s$ of the branch
locus of $f$ is equal to the sum of the ramification indices of points of the
fibres over $s$ of $f$ minus the number of points of the fibre.
\begin{proof}
This is clear for the ramification locus and for the branch locus it follows
from the fact that locally at a closed point branch divisor is the sum of the
norms of the ramification divisors at the points of the fibres and that the norm
of a closed point considered as a divisor equals to the image point which is
seen by looking at valuations of a defining element.
\end{proof}
\end{proposition}
We are now ready to give a description of the deformation theory of Abel maps
(and equivalently Abel curves). To simplify descriptions, for an Abel map
$(\map{f}{P}{Q},C,D)$ over a base $S$ by its \Definition{assigned branch points}
we shall mean the divisor of $Q$ which is the sum of the involutive locus and
the $\infty$-divisor.
%;;%%%Deformation description Theorem Abel deformations
\begin{theorem}\label{Abel deformations}
Let $(\map{f}{P}{Q},C,D)$ be a split Abel map over an algebraically closed field
$\k$. Let $S \subset P(\k)$ be the ramification points that do not map to the
assigned branch points, let $S^e_1$ and $S^e_2$ resp.\ $S^o_1$ and $S^o_2$ be
the ramification points over $(1\co1)$ and $(-1\co1)$ with even resp.\ odd
ramification index (wrt to $f$) and let $T :=
(S,S^o_1,S^e_1,S^o_2,S^e_2)$. Finally, let $r$ associate to a point its
ramification index wrt the map $f$. Then the completion of the local ring of the
stack of Abel maps at the Abel map is isomorphic to the completion of the local
ring of $\Pol{T}$ at the origin and the isomorphism may be assumed to be
stratification preserving.
\begin{proof}
We shall give an isomorphism of deformation functors so we consider a
deformation of the given Abel map over a local Artinian ring $R$ with residue
field $\k$. Note that as $2$ is invertible, the involutive bundle has just the
trivial deformation so we may restrict ourselves to split Abel maps. If we just
consider deformations of the map $f$, then \(curve map deformations) shows that
such deformations are in bijection with tuples $(q_s)_{s \in S'}$, $S'$ being as
in Definition \ref{def:Pol}, where $q_s$ is a normalised polynomial over $R$ of
degree $r(s)$. It remains to understand the influence the choice of relative
Cartier divisors has. Now, the Weierstrass preparation theorem is equivalent to
saying that the ideal of a relative Cartier divisor of $\Sp \pow[R]{x} \to \Sp
R$ is generated by a unique Weierstrass polynomial (i.e., of the form
$x^n+a_1x^{n-1}+\cdots+a_0$ with $a_i \in \maxid_R$) and it is clear from the
uniqueness that inclusion of divisors correspond to divisibility of polynomials
and addition of divisors correspond to product of divisors.

Hence for $s \in S^o_1\disjunion S^e_1\disjunion S^o_2\disjunion S^e_2$ the
inverse image the assigned branch points at $s$ is defined by $q_s$, $C_s$ being
\'etale is defined by a polynomial of degree $0$ or $1$ depending on whether the
degree of $q_s$ is odd or even (as the difference is even). In the even case, if
$p_s$ is the polynomial of $D_s$ we have $q_s=p^2_s$ and $p_s$ is normalised as
$q_s$ is. In the odd case, if $C_s$ is given by $x-a_s$ and $D_s$ by $p_s$ we
have that $q_s=(x-a_s)p^2_s$ and as $q_s$ is normalised we have that $a_s$ is
twice the next to highest coefficient of $p_s$. This shows that the $q_s$ for $s
\in S$ and the $p_s$ for $s$ in the complement gives an $R$-point of $\Pol{T}$
and the converse is also clear. The definition of the stratification of
$\Pol{T}$ has been set up so that the constructed isomorphism preserves the
strata.
\end{proof}
\end{theorem}
\begin{remark}
The assumption of an algebraically closed field as base is just for notational
convenience as is the existence of a splitting.
\end{remark}
We put the most important consequences of this theorem in the following
corollary. Note that we have identified the stack of Abel maps with that of Abel
curves.
%;;%%%Dimension and open strata Corollary Abel smoothness(1) open simple(2)
\begin{corollary}\label{Abel smoothness}
\part $\Ab{g}{n}$ is a smooth stack everywhere of dimension $g$.

\part\label{open simple} The open substack of $\Ab{g}{n}$ consisting of Abel
maps with only simple ramification (i.e., all ramification indices are $\le 2$)
and for which for all branch points outside of the assigned branch points there
is only one ramification point above it is dense.
\begin{proof}
The map $\sAb{g}{n} \to \sAb{g}{n}$ is an \'etale cover so we may deal with the
split case instead. The smoothness follows immediately from the theorem and we
postpone the calculation of the dimension. For the second part we may complete
the local ring at a point and then transfer the problem to the complete local
ring at the origin of $\Pol{T}$. We shall now show that the set of $\Pol{T}$
where the corank of the trace form is $0$ is non-empty and by \(index and
corank) it is enough to show that generically on $\Pol{T}$ the derivative of
each $q_s$ has no multiple roots. If $s \in S$ this is clear as then $q_s$ is a
generic monic polynomial and then so is $1/nq'_s$, where $n$ is the degree of
$q_s$. If $s \in S^e_i$, $i=1,2$, then $q_s=p^2_s$ where $p_s$ is a generic
monic polynomial and thus $q'_s=2p_sp'_s$. Generically $p_s$ and $p'_s$ have no
roots in common, $p_s$ has no double roots and neither has $p'_s$ by the
argument just given. If $s \in S^e_i$ then $q_s=(x-a_s)p^2_s$ where $p_s$ is a
generic monic polynomial with next to highest coefficient $a_s$. Then we have
$q'_s=p_s(p_s+2(x-a_s)p'_s)$. Generically $p_s$ has no double roots. Roots that
are common to $p_s$ and $p_s+2(x-a_s)p'_s$ are also roots of either $x-a_s$ but
generically $a_s$ is not a root of $p_s$ or of $p'_s$ which again is not the
case generically. We are left with showing that generically $p_s+2(x-a_s)p'_s$
has no double roots. Now, $p_s+2(x-a_s)p'_s$ divided by $2n+1$, $n$ being the
degree of $p_s$, is a generic monic polynomial. Indeed, it is easily seen that
the coefficients of $p_s$ can be expressed as polynomials in those of
$p_s+2(x-a_s)p'_s$.

We have thus shown that generically all ramification points are simple and it
remains to show that away from the assigned branch points there is generically
only one ramification point above one branch point. For this we note that for a
given $s \in S$ the contribution from that ramification point to branch locus is
defined by the norm of $q'_s$ and hence what needs to be shown is that for two
$s, s' \in S$ that map to the same point under $f$ the two norms do not have a
common component. Now $q_s$ and $q_{s'}$ are generic polynomials with
independent coefficients. Hence the locus defined by the common components would
have to be independent of both the coefficients of $q_s$ and $q_{s'}$ (and of
course only depend on their union) and would hence have to be constant. At the
origin the full ramification loci consist just of $0$ and so the common locus
would have to be $0$ everywhere. However, $q_s$ has generically no factor in
common with $q'_s$.

Finally, to compute the dimension we may by what has just been proved look only
at the case where all the ramification is simple and outside of the assigned
branch points there is only one ramification point over a given branch point. We
may also assume that the Abel curve is split. Now, if $s \in S$ is a
ramification point then the local deformation at that point depends on one
parameter as $q_s$ is a normalised second degree polynomial whereas for a
ramification point over the assigned branch points the local deformation at that
point depends on zero parameters as $p_s$ is a normalised first degree
polynomial. Hence the dimension is equal to the cardinality of $S$. Let now
$e_i$, $i=1,2$, be the number of ramification points over $\pm1$, let
$e':=e_1+e_2$, and let $e$ be the number of ramification points not above
$\pm1$. By the Hurwitz formula applied to $f$ we have
\begin{displaymath}
-2=-2n+n-1+e'+e
\end{displaymath}
and by the Hurwitz formula applied to double covering ramified at the
non-ramification points over $\pm1$ we have
\begin{displaymath}
2g-2=-4+(2n-2e')
\end{displaymath}
and elimination gives $e=g$.
\end{proof}
\end{corollary}
\end{section}
%;;%%Lyashko-Looijenga map
\begin{section}{The Lyashko-Looijenga map}

By the Lyashko-Looijenga map is generally meant the map that to a family of
finite maps between smooth curves associates the branch locus of each
member. Sometimes one restricts oneself to families where the trace corank of
the branch locus is constant and then it is natural to consider the fibrewise
reduced subscheme of the branch locus. Furthermore, sometimes some of the branch
points are by assumption fixed and then of course it is natural to exclude them
from consideration. Our situation is of this type as the involutive fixed points
are essentially fixed (i.e., they can not move non-trivially in a continuous
fashion) and actually fixed in the split case.

We shall see that the situation is not completely straightforward; our strata on
which the LL-map is defined will generally turn out to be non-reduced which
certainly kills all hope of the LL-map being \'etale. All is not lost however as
a stratum is locally the product of a smooth stack and a zero-dimensional one
and the LL-map turns out to be \'etale on the reduced substack. The most obvious
reason for the stratum being non-reduced is our definition of the branch
locus. This definition is however more or less forced upon us if one wants the
branch divisor to vary continuously (i.e., be a relative Cartier divisor) as
generically the branch divisor is \'etale and hence determined by the condition
that its support be the branch locus.
%;;%%%Ramification specification Definition ramspec
\begin{definition}\label{ramspec}
\part A \Definition{ramification specification of degree $n$} consists of a
finite multiset $S$ of multisets of (strictly) positive integers such that for
each multiset $s$ in $S$ $\sum_{e \in s}e=n$. The \Definition{multiplicity
multiset} associated to $S$ is the multiset \set{\rho(s)}{s \in S}, where
$\rho(s) = \set{e-1}{ e \in S, e > 1}$. The \Definition{total ramification} of
$S$ is 
\begin{displaymath}
\sum_{e \in s \in S}(e-1).
\end{displaymath}

\part An \Definition{Abel ramification specification of order $n$} is a
ramification specification $S$ of order $n$ whose total ramification equals
$n-1$ together with the choice of a submultiset $T$ of $S$ of cardinality
$2$. If $t$ is the number of odd integers, counted with multiplicity, of the
members of $T$, then the \Definition{genus} of $S$ is equal to $(t-4)/2$.
\end{definition}
\begin{remark}
A ramification specification is determined by its associated multiplicity
multiset and the degree $n$. A multiplicity multiset is the same as a
\emph{passport} of \cite{zvonkin99::oen+lyash+looij}.
\end{remark}
We shall now consider stratifications of $\Ab{g}{n}$. First we consider the
ramification stratification of the universal map of $\Ab{g}{n}$ giving a
multiset of multiplicities associated to each stratum. Then we consider its
intersection with the trace stratification associated to the branch locus which
gives a further division of the multiplicities according to which branch point
they are mapped to. This gives exactly an Abel ramification specification $S$ of
order $n$ associated to each such stratum. Conversely, for each Abel
ramification specification $S$ of order $n$ we denote by $\Ab[S]{g}{n}$ the
corresponding stratum.

If $X \to S$ is a map of algebraic stacks and $m$ a positive integer then
$\Conf{X/S}{m}$ (or just $\Conf{X}{m}$ if $S$ is understood) is the
\Definition{$m$-point configuration space}, i.e., the stack quotient by the
permutation action of the symmetric group $\Sigma_m$ on the open substack of the
$m$'th fibre power of $X \to S$ consisting of distinct points.

Our main use of this construction is to the universal involutive projective bundle;
namely the projective bundle $\P \to B\Sigma_2$ that to an involutive line
bundle over $S$ (i.e., a map $S \to B\Sigma_2$) associates the involutive
projective bundle. We then let $\P' \to B\Sigma_2$ be the open substack of $\P$
obtained by removing the section of infinity and the involutive fixed point
set. Note that the universal involutive line bundle of $\Ab{g}{n}$ gives a map
$\Ab{g}{n} \to B\Sigma_2$ and the base of the universal Abel map is just the
pullback of $\P$ under this map.
%;;%%%Def of the LL map
\begin{definition}
Let $S$ be an Abel ramification specification $(S,T)$ of order $n$ and let $m$
be the cardinality of $S$ minus $2$. We define the Lyashko-Looijenga map
\map{LL}{\Ab[S]{g}{n}}{\Conf{\P'/B\Sigma_2}{m}} by associating to an Abel map
over $S$ its reduced branch locus minus assigned base points.
\end{definition}
We can now prove the major result on the LL-map after we have proven the
following lemma.
%;;%%%Odd assigned nilpotency Lemma odd nilpotency
\begin{lemma}\label{odd nilpotency}
Let $R$ be a commutative ring which contains $\Q$ and $a \in R$. Then for a
strictly positive
integer $n$ the polynomial $(t-a^{2n+1})/(t-a)$ is a square of a polynomial
precisely when $a^{n+1}=0$.
\begin{proof}
In the ring of Laurent power series in $t^{-1}$, $R((t^{-1}))$ the polynomial
has the unique square root
\begin{displaymath}
t^{n}\sqrt{1+at^{-1}+\cdots+a^{2n}t^{-2n}}
\end{displaymath}
and hence the polynomial has a polynomial square root precisely when all powers
beyond $t^{-n}$ have zero coefficients in
$\sqrt{1+at^{-1}+\cdots+a^nt^{-2n}}$. This series is obtained by substituting $s
\mapsto at^{-1}$ in $\sqrt{(1-s^{2n+1})/(1-s)}$ which makes it clear that if
$a^{n+1}=0$ then the square root is a polynomial. It is equally clear that the
converse is true if the coefficient of $s^{n+1}$ in $\sqrt{(1-s^{2n+1})/(1-s)}$
is non-zero. However, as $n>0$, $2n+1\ge n+2$, and thus modulo $s^{n+2}$
$\sqrt{(1-s^{2n+1})/(1-s)}$ is congruent to $(1-t)^{-1/2}$ which clearly has all
of its coefficients non-zero.
\end{proof}
\end{lemma}
%;;%%%LL covering map Theorem LL covering
\begin{theorem}\label{LL covering}
Let $S$ be an Abel ramification specification $(S,T)$ of order $n$ and let $m$
be the cardinality of $S$ minus $2$. 

\part[i] The completion of $\Ab[S]{g}{n}$ any geometric point $s=\Sp\k$ is
isomorphic to
\begin{displaymath}
\prod_{2n+1 \in t \in T;\,n >0}\Spf\pow[\k]{a}/(a^{n+1})\times\prod_{s \in
S\setminus
T}\Spf\pow[\k]{\sigma,a_1,\dots,a_{m(s)}}/(\sigma_1,\dots,\sigma_e)
\end{displaymath}
where $m(s) := |\set{e}{e \in s; e \ge 2}|$, $e(s):=\sum_{e \in s}(e-1)$ and 
\begin{displaymath}
\prod_i(s-a_i)^{e_i-1}= s^e +\sum_{1 \le j \le e}(-1)^j\sigma_js^{e-j}
\end{displaymath}
as polynomials in $s$. In particular, $\Ab[S]{g}{n}$ is smooth precisely when
there is exactly one ramification point over each unassigned branch point and no
ramification point of odd ramification index above involutive fix points. It is
always the case that the reduced substack $(\Ab[S]{g}{n})^{red}$ is smooth.

\part[iii] The Lyashko-Looijenga map
\map{LL}{(\Ab[S]{g}{n})^{red}}{\Conf{\P'/B\Sigma_2}{m}} restricted to the
reduced subscheme is an \'etale covering map.
\begin{proof}
We start by making a local calculation. It is clear from Theorem \ref{Abel
deformations} that we get a product over the elements of $S$. Let us first
consider an unassigned branch point. Let $\{e_1,\dots,e_k\}$ be an element of
$S\setminus T$ with the members equal to $1$ removed. Hence a deformation over a
local Artinian ring $R$ is given by a collection $(p_i)_{1\le i \le k}$ of
normalised polynomials with $\deg p_i=e_i$. Now the condition, that the
deformation stay inside the stratum given by $\{e_1,\dots,e_k\}$ means according
to \(Cartier fred) and the identification of Cartier divisors with Weierstrass
polynomials that each $p'_i$ has the form $e_i(x-\alpha_i)^{e_i-1}$ and as $p_i$
is normalised we get that $\alpha_i=0$ and hence that
$p_i(x)=x^{e_i}+b_i$. Furthermore the ramification divisor is defined by
$x^{e-1}$. To compute the branch divisor we have to compute the norm of
$x^{e-1}$ and using the multiplicativity of the norm it is enough to compute the
norm of $x$. Now, it is clear that under the map $\pow[R]{t} \to \pow[R]{x}$
given by $p_i$ we have that $\pow[R]{x}$ is isomorphic to
$\pow[R]{t,x}/(x^{e_i}+b_i-t)$ which gives that the norm of $x$ is
$\pm(t-b_i)$. Hence the branch divisor is given by
$\prod_i(t-b_i)^{e_i-1}$. Now, we are working in the stratum where the fibrewise
reduced branch divisor exists which means that there is a $\sigma \in \maxid_R$
such that $\prod_i(t-b_i)^{e_i-1}=(t-\sigma)^e$, where $e =
\sum_i(e_i-1)$. Comparing next to highest coefficients gives $\sigma =
\sum_i(e_i-1)b_i$ and changing variable $s = t - \sigma$ and putting
$a_i=b_i-\sigma$ gives us $\prod_i(s-b_i)^{e_i-1}=s^e$. This shows that the
universal $R$ is
$\pow[\k]{\sigma,a_1,\dots,a_k}/(\sigma_1,\dots,\sigma_e)$. Now, for degree
reasons as soon as $k > 1$
$\pow[\k]{\sigma,a_1,\dots,a_k}/(\sigma_1,\dots,\sigma_e)$ is strictly larger
than $\pow[\k]{\sigma}$. On the other hand putting $s$ equal to $a_i$ gives
$a_i^e=0$ which shows that when dividing out by the nilradical of
$\pow[\k]{\sigma,a_1,\dots,a_k}/(\sigma_1,\dots,\sigma_e)$ this ring equals
$\pow[\k]{\sigma}$.

Considering now instead one of the involutive fixed points again as we are in a
fixed ramification stratum we get the form $p_i=x^{e_i}+b_i$. This time however
we have that when $e_i$ is even $p_i$ is a square and when it is odd $p_i$ is a
square times a linear polynomial. In the first case it is easy to see that if
$x^{e_i}+b_i$ is a square then $a_i=0$. In the second case, if
$x^{e_i}+b_i=(x-a_i)q^2(x)$ then setting $x=a_i$ we get $b_i=-a_i^{e_i}$ so that
$(x^{e_i}-a_i^{e_i})/(a-a_i)=q^2(x)$ and we conclude from lemma \ref{odd
nilpotency} that this is possible precisely when $a_i^{(e_i+1)/2}=0$. This
concludes the proof of \DHrefpart{i}.

Turning to the \DHrefpart{iii}, that the map is \'etale is clear from the local
calculation as the fibrewise reduced branch divisor is defined by $t-\sigma$,
using the notations of the first part. It remains to prove that it is proper and
for that we shall use the valuative criterion and as everything is of finite
type over $\Sp\Q$ we may restrict ourselves to discrete valuations which we may
assume to be strictly Henselian. Hence we may assume that the map is split and
by \(Pell versality) we may assume that it is given by a Pell family $(P,Q,R)$
such that $P^2-RQ^2=1$ and $P$ is then the Abel map in question. By for instance
\cite[Lemma 3.1]{zvonkin99::oen+lyash+looij} (and the fact also noted in [loc.\
cit.] that the inverse image of the origin under the LL-map is the origin) $P$
has coefficients in $R$. By Gauss lemma so does $R$ and $Q$. The next step is
to show that the discriminant of $R$ is a unit. For this one may reduce modulo
the maximal ideal of $R$ and apply the Hurwitz formula to the map given by
$P$. Indeed, by assumption the number of branch points of $P$ is fixed and hence
by Hurwitz formula the number of ramification points is also fixed. This makes
it impossible for zeros of $R$ to come together.
\end{proof}
\end{theorem}
\begin{remark}
\part The local description of the stratum contradicts \cite[Prop.\
A.3]{ekedahl01::hurwit} which claims that the LL-map always is \'etale. In view
of the theorem (very slightly modified to fit into the context of [loc.\ cit.]) 
this is now seen to be false when there is more than one ramification point over
a branch point. It thus has to be modified to saying that the restriction to the
reduced subscheme of an equisingular stratum is \'etale. Luckily, this is what
is used in the main text and it is also given a topological proof in \cite[Thm.\
4.2]{ekedahl01::hurwit}.

\part It is possible to get a natural interpretation of the reduced structure on
the strata. This will be treated elsewhere.
\end{remark}
\end{section}
%;;%%Topological construction
\begin{section}{Topological construction}

In this section we shall study the covering given by the LL-map. Even though we
can easily get a description for all strata, using the same methods, we shall
only deal with the open stratum as that gives a combinatorial algorithm for
computing the number of connected components of the stack of Abel curves. Note
however that other strata are also interesting. For instance the lowest stratum
where all the branch points are assigned has been considered in connection with
Grothendieck's ``dessin d'enfants'' (cf., \cite{shabat94::plane}).

As usual the fibres of the LL-map are in bijection with conjugacy classes of
certain sequences of the symmetric group. In order to give a procedure for
computing the number of components of a stratum we need generators for the
fundamental group of the appropriate configuration space. The following result
can most certainly be extracted from the literature but for the convenience of
the reader as well as the author we give a proof. We start by giving some
notation. If $C$ is a simple non-closed oriented curve in $\C$ and $S \subset C$
is a finite set then the orientation of $C$ induces a total order on $S$. If $i$
is a positive integer strictly smaller than the number $s$ of elements of $S$
then we define, as usual, the elements $\sigma_i$ of the braid group on $s$
strands given by letting the $i$'th point move along $C$ to the position of the
$i+1$'st point to the right of $C$ and letting the $i+1$'st point move along $C$
to the position of the $i$'th point to the left of $C$ (``right'' and ``left''
being from the point of view of the orientation of $C$). If $s \ge 3$ and we
define $\tau_1$ resp.\ $\tau_2$ to be the braids that takes the second resp.\
$s-1$'st point and moves along $C$ on the left resp.\ right hand side till just
before the first resp.\ last point, then circles that point once
counter-clockwise, and returns back to its original position along the right
resp.\ left hand side of $C$. (They are equal $\sigma_1^2$ resp.\
$\sigma_{s-1}^{-2}$.)
%;;%%%Braid generators Proposition braid generators
\begin{proposition}\label{braid generators}
Let $C:=[-1,1]$ oriented in any direction and $S \subset C$ a finite subset with
$s$ elements containing $\pm1$ and let $A := \C\setminus\{\pm1\}$. Then the map
induced by the inclusion $\{x_1,\dots,x_{s-2}\} \mapsto
\{-1,1,x_1,\dots,x_{s-2}\}$
\begin{displaymath}
\pi_1(\Conf{A}{s-2},S\setminus\{\pm1\}) \to \pi_1(\Conf{\C}{s},S)
\end{displaymath}
is an injection whose image is generated by $\sigma_i$, $2 \le i \le s-2$,
$\tau_1$ and $\tau_2$.
\begin{proof}
By possibly applying $z \to -z$ we may assume that the orientation of $[-1,1]$
is such that $-1$ becomes its first element. 

To begin with it is clear the $\sigma_i$, $2 \le i \le s-2$, $\tau_1$ and
$\tau_2$ lie in the image. Recall that we have a surjection
$\pi_1(\Conf{\C}{s},S) \to \Sigma_s$ taking $\sigma_i$ to the transposition
$(i-1,i)$. The image of $\pi_1(\Conf{s-2}{A},S\setminus\{\pm1\})$ maps into
$\Sigma_{s-2}$ considered as the subgroup that fixes the first and last elements
and as $\sigma_i$ maps to $(i-1,i)$ the subgroup generated by them maps
surjectively onto $\Sigma_{s-2}$. Hence both for injectivity and generation it
suffices to consider $\pi_1(\oConf{A}{s-2},S\setminus\{\pm1\}) \to
\pi_1(\oConf{\C}{s},S)$, where $\oConf{X}{t}$ is the space of ordered
$t$-subsets of $X$ and to show that the map is injective and the image is
generated by the conjugates of $\sigma_i^2$, $1 \le s \le s-1$ in the group
generated by $2 \le i \le s-2$, $\tau_1$ and $\tau_2$.

Now, by conjugating by the $\sigma_i$, $2 \le i \le s-1$ we can get from
$\tau_1$ and $\tau_2$ all braids $A^1_i$ and $A^2_i$ defined like $\tau_1$
resp.\ $\tau_2$ only starting at the $i$'th point for $2 \le i \le s-1$ as well
as the $A_{ij}$, $2 \le i < j \le s-2$ defined like $\tau_1$ only starting at
the $j$'th point and encircling the $i$'th point. (The $A_{ij}$ are the $A_{ij}$
of \cite[1-11]{birman74::braid}, $A^1_i$ is $A_{1i}$ and $A^2_i$ is a mirror
image of $A_{s-i,s}$.) Our aim is to show the injectivity and that these
elements generate the image. In this we shall follow the proof of \cite[Lemma
1.8.2]{birman74::braid} and we start following \cite{birman74::braid} in using
the notation $F_{m,n}(X)$ for $\oConf{X\setminus Q_m}{n}$ where $Q_m$ is a fixed
subset of $X$ of cardinality $m$ and will use of the theorem of Fadell and
Neuwirth (cf., \cite[Thm.\ 1.2]{birman74::braid}) which says that when $X$ is a
manifold, then the projection on the first $r$ factors $F_{m,n} \to F_{m,r}$ is
a fibration with fibre $F_{m+r,n-r}$. Applied to $r=n-1$ and $X=\C$ and $X=A$
this will allow us to prove the statement by induction. As the involved spaces
are acyclic, the fibrations give short exact sequences and by induction we are
reduced to showing that for $1 \le i < s$
\begin{displaymath}
\pi_1(A\setminus S_i,x_{i+1}) \to \pi_1(\C\setminus S_i,x_{i+1}),
\end{displaymath}
where $S_i$ consists of the $i$ first elements of $S$ and $x_{i+1}$ is the
$i+1$'st element,
is an injection and that the image is contained in the subgroup generated by 
$A_{k,i+1}$ and $A^1_{i+1}$ (and when $i=s-1$ also the $A^2_k$) for $1 \le i \le
k$. This however is clear.
\end{proof}
\end{proposition}
This result combined with Theorem \ref{LL covering} and Corollary \ref{open
simple} allows us to give a combinatorial description of the number of
components of $\Ab{g}{n}$ and $\sAb{g}{n}$. For this we first introduce the
following definition.
%;;%%%Def of combinatorial data
\begin{definition}
Let $N_{g,n}$ be the set of tuples $(\sigma,\sigma_1,\dots,\sigma_g,\tau) \in
(\Sigma_n)^{g+2}$ fulfilling the conditions
\begin{itemize}
\item $\sigma\sigma_1,\dots,\sigma_g\tau$ is an $n$-cycle and

\item the $\sigma_i$ are transpositions and $\sigma$ and $\tau$ are products of
disjoint transpositions and the sum of the number of fixed points of $\sigma$
and of $\tau$ equals $2g+2$.
\end{itemize}
Let $M_{g,n}$ be the set of orbits of the action of $\Sigma_n$ on $N_{g,n}$
given by 
\begin{displaymath}
(\rho,(\sigma,\sigma_1,\dots,\sigma_g,\tau)) \mapsto
(\rho\sigma\rho^{-1},\rho\sigma_1\rho^{-1},\dots,
\rho\sigma_g\rho^{-1},\rho\tau\rho^{-1}).
\end{displaymath}
\end{definition}
Thus armed we can give a combinatorial description of the set of connected
components of the stacks of (split) Abel curves.
%;;%%%Number of components
\begin{theorem}
The set of connected components of $\sAb{g}{n}$ is in bijection with
equivalence classes of $M_{g,n}$ under the equivalence relation generated by the
relations.
\begin{itemize}
\item $(\sigma,\sigma_1,\dots,\sigma_i,\sigma_{i+1},\dots,\sigma_g,\tau) \sim
(\sigma,\sigma_1,\dots,\sigma_i\sigma_{i+1}\sigma_i^{-1},\sigma_i,\dots,\sigma_g,\tau)$
for all $1 \le i < g$.

\item $(\sigma,\sigma_1,\dots,\sigma_g,\tau) \sim
(\sigma[\sigma_1,\sigma],\sigma\sigma_1\sigma^{-1},\dots,\sigma_g,\tau)$, where
$[\sigma_1,\sigma]=\sigma_1\sigma\sigma_1^{-1}\sigma^{-1}$.

\item $(\sigma,\sigma_1,\dots,\sigma_g,\tau) \sim
(\sigma,\sigma_1,\dots,\tau^{-1}\sigma_g\tau,[\tau^{-1},\sigma_g^{-1}]\tau)$.
\end{itemize}
The set of connected components of $\Ab{g}{n}$ is in bijection with
equivalence classes of $M_{g,n}$ under the equivalence relation generated by the
above relations together with the relation
\begin{displaymath}
(a_1,a_2,\dots,a_{g+2}) \sim (b_{g+2}a_{g+2}b_{g+2}^{-1},\dots,b_2a_2b_2^{-1},a_1),
\end{displaymath}
where $b_i=a_1\dots a_{i-1}$ for $i \ge 2$.
\begin{proof}
The part on $\sAb{g}{n}$ follows directly from the fact that the LL-map is an
\'etale covering (Theorem \ref{LL covering}), that the fibres of the LL-mapping
are in bijection with $M_{g,n}$, the description of the generators for the
fundamental group for the target of the LL-map (Proposition \ref{braid
generators}) and the formula for the action of the $\sigma_i$ on $M_{g,n}$.

As for the $\Ab{g}{n}$-part the LL-map has as target the quotient of
$\Conf{\A^1\setminus\{\pm1\}}{g}$ divided by the map induced by $z \mapsto
-z$. Hence we have to add the relation that identifies an equivalence class of
maps from the fundamnetal group to $\Sigma_n$ with the one obtained by composing
with the action of the (outer) automorphism induced by $z \mapsto -z$. For that
we choose as basepoint of $\Conf{\A^1\setminus\{\pm1\}}{g}$ the set
$\{-1/2,-1/3,\dots,1/3,1/2\}$ (with $0$ included if $g$ is odd) and as basepoint
for $\A^1\setminus\{-1,-1/2,-1/3,\dots,1/3,1/2,1\}$ $i$. Acting by $z \mapsto
-z$ gives us $-i$ as new basepoint and we identify fundamental groups by
choosing a curve from $-i$ to $i$ going to the left of $\{-1,-1/2,-1/3,\dots,1/3,1/2,1\}$.
\end{proof}
\end{theorem}
\end{section}
\bibliography{preamble,abbrevs,alggeom,algtop,ekedahl}

\newcommand\eprint[1]{\texttt{arXiv:#1}}\def\cprime{$'$}
\providecommand{\bysame}{\leavevmode\hbox to3em{\hrulefill}\thinspace}
\begin{thebibliography}{ELSV01}

\bibitem[Ab26]{abel26::sur+r+r}
N.~H. Abel, \emph{Sur l'int{\'e}gration de la formule diff{\'e}rentielle
  $\frac{\rho\,dx}{\sqrt r}$, $r$ et $\rho$ {\'e}tant des fonctions
  ent{\`e}res.}, J. Reine Angew. Math. (1826), no.~1, 105--144.

\bibitem[Bi74]{birman74::braid}
J.~S. Birman, \emph{Braids, links, and mapping class groups}, Princeton
  University Press, Princeton, N.J., 1974, Annals of Mathematics Studies, No.
  82.

\bibitem[ELSV01]{ekedahl01::hurwit}
T.~Ekedahl, S.~Lando, M.~Shapiro, and A.~Vainshtein, \emph{Hurwitz numbers and
  intersections on moduli spaces of curves}, Invent. Math. \textbf{146} (2001),
  no.~2, 297--327, \eprint{math.AG/0004096}.

\bibitem[SZ94]{shabat94::plane}
G.~Shabat and A.~Zvonkin, \emph{Plane trees and algebraic numbers}, Jerusalem
  combinatorics '93, Contemp. Math., vol. 178, Amer. Math. Soc., Providence,
  RI, 1994, pp.~233--275.

\bibitem[ZL99]{zvonkin99::oen+lyash+looij}
D.~Zvonkin and S.~K. Lando, \emph{On multiplicities of the
  {L}yashko-{L}ooijenga mapping on strata of the discriminant}, Funktsional.
  Anal. i Prilozhen. \textbf{33} (1999), no.~3, 21--34.

\end{thebibliography}
\bibliographystyle{pretex}

\end{document}